\documentclass{mfat}
\pagespan{184}{200}

\usepackage{mmap}
\usepackage[utf8]{inputenc}

\theoremstyle{plain}
\newtheorem{theorem}{Theorem}
\newtheorem{lemma}[theorem]{Lemma}
\newtheorem{proposition}[theorem]{Proposition}
\newtheorem{statement}[theorem]{Statement}
\newtheorem{corollary}[theorem]{Corollary}
\theoremstyle{definition}

\newtheorem{remark}[theorem]{Remark}

\theoremstyle{plain}
\newtheorem*{theorem*}{Theorem}
\newtheorem*{theorem*A}{Theorem A}
\newtheorem*{theorem*B}{Theorem B}
\newtheorem*{theorem*C}{Theorem C}
\newtheorem*{lemma*}{Lemma}
\newtheorem*{proposition*}{Proposition}
\newtheorem*{statement*}{Statement}
\newtheorem*{corollary*}{Corollary}
\theoremstyle{definition}
\newtheorem*{definition*}{Definition}
\theoremstyle{remark}
\newtheorem*{notation*}{Notation}
\newtheorem*{remark*}{Remark}

\begin{document}

\title[1-D Schr\"{o}dinger operators with singular periodic potentials]
      {One-dimensional Schr\"{o}dinger operators with singular periodic potentials}

\author{Vladimir Mikhailets}
\address{Institute of Mathematics, National Academy of Sciences of Ukraine, 3 Tereschenkivs'ka, Kyiv, 01601, Ukraine}
\email{mikhailets@imath.kiev.ua}

\author{Volodymyr Molyboga}
\address{Institute of Mathematics, National Academy of Sciences of Ukraine, 3 Tereschenkivs'ka, Kyiv, 01601, Ukraine}
\email{molyboga@imath.kiev.ua}

\subjclass[2000]{Primary 34L05, 47A05; Secondary 34L40, 47A10, 47B25}
\date{13/03/2008}
\dedicatory{Dedicated to M. L. Gorbachuk on the occasion of his
70th birthday.}
\keywords{Hill equations, Schr\"{o}dinger operators, singular
potentials, spectrum gaps, periodic eigenvalues.}

\begin{abstract}
  We study the one-dimensional Schr\"{o}dinger operators
  \begin{equation*}
    S(q)u:=-u''+q(x)u,\quad u\in \mathrm{Dom}\left(S(q)\right),
  \end{equation*} 
  with $1$-periodic real-valued singular potentials $q(x)\in
  H_{\operatorname{per}}^{-1}(\mathbb{R},\mathbb{R})$ on the Hilbert
  space $L_{2}\left(\mathbb{R}\right)$. We show equivalence of five
  basic definitions of the operators $S(q)$ and prove that they are
  self-adjoint. A new proof of continuity of the spectrum of the
  operators $S(q)$ is found. Endpoints of spectrum gaps are
  precisely described.
\end{abstract}

\maketitle

\section{Introduction}\label{sec_Int}
On the complex Hilbert space $L_{2}(\mathbb{R})$, we consider the
one-dimensional Schr\"{o}dinger operators
\begin{equation}\label{eq_10_Int}
  S(q)u:=-u''+q(x)u,\quad u\in \mathrm{Dom}\,\left(S(q)\right),
\end{equation}
with real-valued $1$-periodic distribution potentials $q(x)$, the so-called the
Hill-Schr\"{o}dinger operators.

Assuming that
\begin{equation}\label{eq_12_Int}
  q(x)=\sum_{k\in 2\mathbb{Z}}\widehat{q}(k)e^{i k\pi x}\in
  H_{\operatorname{per}}^{-1}(\mathbb{R},\mathbb{R}),
\end{equation}
that is,
\begin{equation*}
    \sum_{k\in
    2\mathbb{Z}}(1+|k|)^{-2}|\widehat{q}(k)|^{2}<\infty\quad\text{and}\quad
    \widehat{q}(k)=\overline{\widehat{q}(-k)}\quad \forall k\in
    2\mathbb{Z},
\end{equation*}
the Hill-Schr\"{o}dinger operators $S(q)$ can be well defined on the Hilbert
space $L_{2}(\mathbb{R})$ in the following different ways:
\begin{enumerate}
  \item  as minimal/maximal quasi-differential operators
  $S_{\operatorname{min}}(q)$/$S_{\operatorname{max}}(q)$;
  \item  as Friedrichs extensions $S_{F}(q)$ of
  quasi-differential operators $S_{\operatorname{min}}(q)$;
  \item  as form-sum operators $S_{\operatorname{form}}(q)$;
  \item as the limit $S_{\operatorname{lim}}(q)$ of sequences of the
    Hill-Schr\"{o}dinger operators with smooth periodic potentials
    in the norm resolvent sense.
\end{enumerate}

Hryniv and Mykytyuk \cite{HrMk}, Djakov and Mityagin \cite{DjMt}
studied the Friedrichs extensions $S_{F}(q)$, but Korotyaev \cite{Krt}
treated the form-sum operators $S_{\operatorname{form}}(q)$. We propose to
join together these results showing an equivalence of all
definitions.

More precisely, we will prove the following statements.
\begin{theorem*A}{\rm{(Theorem \ref{th_26_MnRs})}}.
The Hill-Schr\"{o}dinger quasi-differential operators $S_{\operatorname{max}}(q)$ with
distributional potentials $q(x)\in
H_{\operatorname{per}}^{-1}\left(\mathbb{R},\mathbb{R}\right)$ are self-adjoint.
\end{theorem*A}

\begin{theorem*B}{\rm{(Corollary \ref{cr_28_MnRs}, Corollary \ref{cr_29_MnRs},
      Theorem \ref{th_34_MnRs})}}. The quasi-differential operators
  $S_{\operatorname{min}}(q)$ and $S_{\operatorname{max}}(q)$,
  the Friedrichs extensions $S_{F}(q)$, the form-sum operators
  $S_{\operatorname{form}}(q)$,  and the ope\-rators
  $S_{\operatorname{lim}}(q)$ coincide.
\end{theorem*B}

In the paper \cite[Theorem 3.5]{HrMk} the authors tried to show that
the operators $S_{\operatorname{max}}(q)$ and $S_{F}(q)$ coincide.
But the proof of this assertion was erroneous. Our proofs of Theorem
A and Theorem B are based on a different idea (see Lemma
\ref{lm_18_Prl}).

The equality $S(q)=S_{\operatorname{lim}}(q)$, together with the
classical Birkhoff-Lyapunov theorem, allow to prove the following
statement.
\begin{theorem*C}{\rm{(Theorem \ref{th_36_MnRs})}}.
  The Hill-Schr\"{o}dinger operators $S(q)$ with distributional
  potentials $q(x)\in
  H_{\operatorname{per}}^{-1}(\mathbb{R},\mathbb{R})$ have
  continuous spectra with the band and the gap structures being such
  that the endpoints $\{\lambda_{0}(q),
  \lambda_{k}^{\pm}(q)\}_{k=1}^{\infty}$ of the spectrum gaps
  satisfy the inequalities
\begin{equation*}
  -\infty<\lambda_{0}(q)<\lambda_{1}^{-}(q)\leq\lambda_{1}^{+}(q)<
  \lambda_{2}^{-}(q)\leq\lambda_{2}^{+}(q)<\cdots\,
\end{equation*}
Moreover, endpoints of the spectrum gaps for even (odd) numbers
$k\in \mathbb{Z}_{+}$ are periodic (semiperiodic) eigenvalues of the
following problem on the interval $[0,1]$:
\begin{equation*}
   S_{\pm}(q)u=-u''+q(x)u=\lambda u,\quad u\in
   \mathrm{Dom}\left(S_{\pm}(q)\right).
\end{equation*}
\end{theorem*C}
It is interesting to remark that the last assertion is nontrivial, and for
more singular $\delta'$-interactions, that is if
\begin{equation*}
    q(x)=\sum_{k\in \mathbb{Z}}\beta\,\delta'(x-k)\notin
    H_{\operatorname{per}}^{-1}(\mathbb{R}),\quad \beta<0,
\end{equation*}
it could still occur that endpoints of the spectrum gaps for even
(odd) numbers $k\in \mathbb{Z}_{+}$ are \textit{semiperiodic
  (periodic)} eigenvalues of the problem on the interval
$[0,1]$, see~\cite[Theorem~III.3.6]{Alb}.

In the closely related paper of Hryniv and Mykytyuk \cite{HrMk}, the
authors have established that spectra of the operators $S(q)$ are
absolutely continuous.

\section{Preliminaries}\label{sec_Prl}
\subsection{Sobolev spaces}\label{ssec_SbSp_Prl}
Let us denote by $\mathfrak{D}'_{1}\left(\mathbb{R}\right)$ the
Schwartz space of $1$-periodic distributions defined on the whole
real axis $\mathbb{R}$ (see \cite{Vld}). To have a detailed
characterization of $1$-periodic distributions, we will use Sobolev
spaces.

Consider the Sobolev spaces
$H_{\operatorname{per}}^{s}(\mathbb{R})$, $s\in \mathbb{R}$, of
$1$-periodic functions (distributions) defined by means of their
Fourier coefficients,
\begin{align*}
  H_{\operatorname{per}}^{s}(\mathbb{R}) & :=\bigg\{f=\sum_{k\in 2\mathbb{Z}}\widehat{f}(k)e^{i k\pi
  x}\left|\;\parallel f\parallel_{H_{\operatorname{per}}^{s}(\mathbb{R})}<\infty\right.\bigg\}, \\
  \parallel f\parallel_{H_{\operatorname{per}}^{s}(\mathbb{R})}&:=\bigg(\sum_{k\in 2\mathbb{Z}}
  \langle k\rangle^{2s}|\widehat{f}(k)|
  ^{2}\bigg)^{1/2},\quad \langle k\rangle:=1+|k|, \\
  \widehat{f}(k)&:=\langle f,e^{i  k\pi x}\rangle_{L_{2,{\operatorname{per}}}(\mathbb{R})}, \quad
   k\in 2\mathbb{Z}, \\
  2\mathbb{Z}&:=\left\{ k\in\mathbb{Z}\; \left|\;k\equiv 0\;
  (\mathrm{mod}\,2)\right.\right\}.
\end{align*}
The sesquilinear form $\langle
\cdot,\cdot\rangle_{L_{2,{\operatorname{per}}}(\mathbb{R})}$ pairs
the dual, respectively $L_{2,{\operatorname{per}}}(\mathbb{R})$,
spaces $H_{\operatorname{per}}^{s}(\mathbb{R})$ and
$H_{\operatorname{per}}^{-s}(\mathbb{R})$, and is an extension by
continuity of the $L_{2,{\operatorname{per}}}(\mathbb{R})$-inner
product \cite{Brz, Grb},
\begin{equation*}
  \langle
  f,g\rangle_{L_{2,{\operatorname{per}}}(\mathbb{R})}:=\int_{0}^{1}f(x)\overline{g(x)}\,dx\quad
  \forall f,g\in L_{2,{\operatorname{per}}}(\mathbb{R}).
\end{equation*}
It should be noted that
\begin{equation*}
  H_{\operatorname{per}}^{0}(\mathbb{R})=L_{2,{\operatorname{per}}}(\mathbb{R}),
\end{equation*}
and we denote by
$\mathfrak{D}'_{1}\left(\mathbb{R},\mathbb{R}\right)$ and
$H_{\operatorname{per}}^{s}(\mathbb{R},\mathbb{R})$, $s\in
\mathbb{R}$, the \textit{real-valued} $1$-periodic
distributions from the correspondent spaces,
\begin{align*}
  \mathfrak{D}'_{1}\left(\mathbb{R},\mathbb{R}\right)&:=\left\{f(x)\in
  \mathfrak{D}'_{1}(\mathbb{R})\left|\,\mathrm{Im}f(x)=0
  \right.\right\}, \\
  H_{\operatorname{per}}^{s}(\mathbb{R},\mathbb{R})&:= \left\{f(x)\in
  H_{\operatorname{per}}^{s}(\mathbb{R})\left|\,\mathrm{Im}f(x)=0
  \right.\right\}.
\end{align*}
Note that $\mathrm{Im}f(x)=0$ for a 1-periodic distribution $f(x)\in
\mathfrak{D}'_{1}\left(\mathbb{R}\right)$ means that
\begin{equation*}
    \widehat{f}(2k)=\overline{\widehat{f}(-2k)}\quad \forall k\in
    \mathbb{Z}.
\end{equation*}

\subsection{Quasi-differential equations}\label{ssec_QDE_Prl}
The differential expressions in the right-hand of \eqref{eq_10_Int},
by introducing quasi-derivatives
\begin{equation*}
  u^{[1]}(x):=u'(x)-Q(x) u(x),
\end{equation*}
can be re-written as quasi-differential expressions \cite{SvSh1,
  SvSh2},
\begin{equation*}
  l_{Q}[u]:=-(u'-Qu)'-Q(u'-Qu)-Q^{2}u,
\end{equation*}
which are well defined if $u, u^{[1]}\in W_{1,{\operatorname {loc}}}^{1}(\mathbb{R})$ \cite{Nai}.
\begin{proposition}{\rm{(Existence and Uniqueness
      Theorem)}}\label{pr_12_Prl}.  Let $\lambda\in \mathbb{C}$ and
  $f(x)\in L_{1,{\operatorname {loc}}}(\mathbb{R})$.  Then, for any
  complex numbers $c_{0},\; c_{1}\in \mathbb{C}$ and arbitrary
  $x_{0}\in \mathbb{R}$, the quasi-differential equation
\begin{equation}\label{eq_20_Prl}
  l_{Q}[u]=\lambda u + f,\quad \lambda\in \mathbb{C},\quad f\in L_{1,{\operatorname {loc}}}(\mathbb{R}),
\end{equation}
has one and only one solution $u\in W_{1,{\operatorname
    {loc}}}^{1}(R)$ satisfying the initial conditions
\begin{equation*}\label{eq_22_Prl}
  u(x)\left|\right._{x=x_{0}}=c_{0},\quad
  u^{[1]}(x)\left|\right._{x=x_{0}}=c_{1}.
\end{equation*}
\end{proposition}

For the quasi-differential equation \eqref{eq_20_Prl} there is a relating
normal $2$-dimensional system of the first order differential
equations with locally integrable coefficients,
\begin{equation*}\label{eq_24_Prl}
  \begin{pmatrix}
    u_{1} \\
    u_{2} \
  \end{pmatrix}
  '=
  \begin{pmatrix}
    Q & 1 \\
    -\lambda-Q^{2} & -Q \
  \end{pmatrix}
  \begin{pmatrix}
    u_{1} \\
    u_{2} \
  \end{pmatrix}
  +
  \begin{pmatrix}
    0 \\
    -f \
  \end{pmatrix}
  ,
\end{equation*}
where $u_{1}(x):=u(x)$, $u_{2}(x):=u^{[1]}(x)$.

Then Proposition \ref{pr_12_Prl} follows from \cite[Theorem 1, \S
16]{Nai}, also see \cite{AhGl}.

\enlargethispage{-1cm}

\begin{lemma}{\rm{(Lagrange formula)}}\label{lm_14_Prl}.
  Let $u(x)$ and $v(x)$ be functions such that the
  quasi-differential expressions $l_{Q}[\cdot]$ are well defined.
  Then the following Lagrange formula holds:
\begin{equation*}\label{eq_26_Prl}
 l_{Q}[u]\overline{v}-u l_{Q}[\overline{v}]=\frac{d}{dx}[u,v]_{x},
\end{equation*}
where the sesquilinear forms $[u,v]_{x}$ are defined by
\begin{equation*}
  [u,v]_{x}:=u(x)\overline{\left(v'(x)-Q(x)v(x)\right)}-\left(u'(x)-Q(x)u(x)\right)\overline{v(x)}.
\end{equation*}
\end{lemma}

\begin{proof}
  It follows at once that $u(x)$ and $v(x)$ are such that
\begin{equation*}
  u,u'-Qu\in W_{1,{\operatorname {loc}}}^{1}(\mathbb{R})\quad\text{and}\quad v,v'-Qv\in
  W_{1,{\operatorname {loc}}}^{1}(\mathbb{R}).
\end{equation*}
Then we have
\begin{align*}
  \frac{d}{dx}[u,v]_{x} & \equiv \frac{d}{dx}
  \left(u\overline{\left(v'-Qv\right)}-\left(u'-Qu\right)\overline{v}\right) \\
  & =u'\overline{\left(v'-Qv\right)}+ u\overline{\left(v'-Qv\right)'}
   -\left(u'-Qu\right)'\overline{v} -\left(u'-Qu\right)\overline{v'} \\
  & =l_{Q}[u]\overline{v}-u l_{Q}[\overline{v}]+Qu'\overline{v}
  -Qu\overline{v'}+u'\overline{\left(v'-Qv\right)}
  -\left(u'-Qu\right)\overline{v'} \\
  & =l_{Q}[u]\overline{v}-u l_{Q}[\overline{v}],
\end{align*}
since it follows from the assumptions that
\begin{equation*}
  u'\overline{v'},\, Q^{2}u\overline{v},\, Qu'\overline{v},\,
  Qu\overline{v'}\in L_{1,{\operatorname {loc}}}(\mathbb{R}).
\end{equation*}

The proof is complete.
\end{proof}

Integrating both sides of the Lagrange formula over the compact
interval $[\alpha,\beta]\Subset \mathbb{R}$ we obtain the Lagrange
identity in an integral form,
\begin{equation}\label{eq_28_Prl}
  \int_{\alpha}^{\beta}l_{Q}[u]\overline{v}\,dx-\int_{\alpha}^{\beta}
  u l_{Q}[\overline{v}]\,dx=[u,v]_{\alpha}^{\beta},
\end{equation}
where
\begin{equation*}
  [u,v]_{\alpha}^{\beta}:=[u,v]_{\beta}-[u,v]_{\alpha}.
\end{equation*}

\subsection{Quasi-differential operators on a finite
interval}\label{ssec_QDOfi_Prl} Here, following Savchuk and Shkalikov
\cite{SvSh1}, we give a brief review of results related to Sturm-Liouville
operators with distribution potentials defined on a finite interval.

On the Hilbert space $L_{2}(0,1)$, we consider the Sturm-Liouville
operators
\begin{equation*}
  L(q)u:=-u''+q(x)u,\quad u\in \mathrm{Dom}\,\left(L(q)\right),
\end{equation*}
with real-valued  distribution potentials $q(x)\in H^{-1}\left([0,1],
\mathbb{R}\right)$, i.e.,
\begin{equation*}
  Q(x)=\int q(\xi)\,d\xi\in L_{2}\left((0,1),\mathbb{R}\right).
\end{equation*}

Set
\begin{align*}
  L_{\operatorname{max}}(q)u & :=l_{Q}[u], \\
  \mathrm{Dom}(L_{\operatorname{max}}(q)) & :=\left\{u\in L_{2}(0,1)\left|\right.u,u'-Qu\in W_{1}^{1}[0,1],\,
  l_{Q}[u]\in L_{2}(0,1)\right\},
\end{align*}
and
\begin{align*}
  \dot{L}_{\operatorname{min}}(q)u & :=l_{Q}[u], \\
  \mathrm{Dom}(\dot{L}_{\operatorname{min}}(q)) & :=\left\{u\in \mathrm{Dom}(L_{\operatorname{max}}(q))\left|\right.
  \mathrm{supp}\,u\Subset [0,1]\right\}.
\end{align*}
We also consider the operators
\begin{align*}
  L_{\operatorname{min}}(q)u & :=l_{Q}[u], \\
  \mathrm{Dom}(L_{\operatorname{min}}(q)) & :=\left\{u\in \mathrm{Dom}(L_{\operatorname{max}}(q))\left|\right.
  u^{[j]}(0)=u^{[j]}(1)=0,\; j=0,1\right\}.
\end{align*}

\begin{proposition}{\rm{(\cite{SvSh1})}}\label{pr_15_Prl}.
   Suppose that $q(x)\in H^{-1}\left([0,1], \mathbb{R}\right)$.
  Then the following statements are true:
\begin{itemize}
  \item
    [(\rm I)] The operators $L_{\operatorname{min}}(q)$ are densely
    defined on the Hilbert space $L_{2}(0,1)$.
  \item
    [(\rm{II})] The operators $L_{\operatorname{min}}(q)$ and
    $L_{\operatorname{max}}(q)$ are mutually adjoint,
  \begin{equation*}
     L_{\operatorname{min}}^{\ast}(q)=L_{\operatorname{max}}(q),\qquad L_{\operatorname{max}}^{\ast}(q)=
     L_{\operatorname{min}}(q).
  \end{equation*}
  In particular, the operators $L_{\operatorname{min}}(q)$ and
  $L_{\operatorname{max}}(q)$ are closed.
\end{itemize}
\end{proposition}

In Statement \ref{st_15.1_Prl}, which is proved in Appendix A.1, we
establish relationships between the operators
$\dot{L}_{\operatorname{min}}(q)$ and $L_{\operatorname{min}}(q)$.
\begin{statement}\label{st_15.1_Prl}
  The operators $L_{\operatorname{min}}(q)$ are closures of the
  operators $\dot{L}_{\operatorname{min}}(q)$,
\begin{equation*}
  L_{\operatorname{min}}(q)=(\dot{L}_{\operatorname{min}}(q))^{\sim}=\dot{L}_{\operatorname{min}}^{\ast\ast}(q).
\end{equation*}
\end{statement}
\section{Main results}\label{sec_MnRs}
\subsection{A principal lemma}
The following operator-theory result is an essential part of our
approach. In this section, we will give two important applications.

\begin{lemma}\label{lm_18_Prl}
  Let $A$ be a linear operator that is densely defined and closed on
  a complex Banach space $X$, and let $B$ be a linear operator
  bounded on $X$ such that
\begin{itemize}
  \item [(a)] $BA\subset AB$ ($A$ and $B$ commute);
  \item [(b)] $\sigma_{p}(B)=\varnothing$ (the point spectrum
    $\sigma_{p}(B)$ of the operator $B$ is empty).
\end{itemize}
Then the operator $A$ has no eigenvalues of finite multiplicity.
\end{lemma}
\begin{proof}
  Suppose that the operator $A$ has an eigenvalue $\lambda\in
  \sigma_{p}(A)$ of finite multiplicity, and let $G_{\lambda}$ be
  the corresponding eigenspace.
  
  Further, let $f$ be an eigenvector of the operator $A$,
\begin{equation*}
  Af=\lambda f,\quad f\in G_{\lambda}.
\end{equation*}
Then
\begin{equation*}
  A(Bf)=B(Af)=\lambda(Bf),\quad f\in
  G_{\lambda},
\end{equation*}
whence we conclude that
\begin{equation*}
  B G_{\lambda}\subset G_{\lambda}.
\end{equation*}
The assumption $\mathrm{dim}(G_{\lambda})\in \mathbb{N}$ implies
that the point spectrum $\sigma_{p}(B)$ of the operator $B$ is not
empty.  This contradicts condition~(b).

The proof is complete.
\end{proof}

\begin{remark}\label{rm_21_MnRs}
  The condition (b) is satisfied if
  $X=L_{p}(\mathbb{R},\mathbb{C})$, $1\leq p<\infty$, and $B$ is a
  shift operator,
\begin{equation*}
  B:\,y(x)\mapsto y(x+T),\quad T>0.
\end{equation*}

Indeed, the operator $B$ is unitary on the space
$X=L_{p}(\mathbb{R},\mathbb{C})$. Therefore,
\begin{equation*}
  \sigma_{p}(B)\subset \sigma(B)=\left\{\lambda\in
  \mathbb{C}\left|\right.|\lambda|=1\right\},
\end{equation*}
and the identity
\begin{equation*}
  B y(x)=\lambda y(x)=y(x+T),\quad y(x)\not\equiv 0,\quad
  |\lambda|=1,
\end{equation*}
implies that the function $|y(x)|$ is $T$-periodic. Then
$y(x)\not\in L_{p}(\mathbb{R},\mathbb{C})$, and we conclude that
$\sigma_{p}(B)=\varnothing$.

Condition (a) means in this case that the operator $A$ is
$T$-periodic on the line.
\end{remark}
\subsection{Self-adjointness of the Hill-Schr\"{o}dinger operators with
  distribution potentials} 

If assumption \eqref{eq_12_Int} is true, then the distribution
potentials $q(x)$ can be represented as
\begin{equation*}
  q(x)=C+Q'(x)
\end{equation*}
with
\begin{equation*}
  C=\widehat{q}(0)
\end{equation*}
and
\begin{equation*}
  Q(x)=\sum_{k\in 2\mathbb{Z}\setminus\{0\}}
  \frac{1}{i k\pi}\widehat{q}(2k)e^{i k\pi x}\in
  L_{2,{\operatorname{per}}}(\mathbb{R},\mathbb{R})
\end{equation*}
such that
\begin{equation*}
  \langle q,\varphi\rangle=-\langle Q,\varphi'\rangle\quad \forall\varphi\in
  C_{\operatorname {comp}}^{\infty}(\mathbb{R}),
\end{equation*}
see \cite[Proposition 1]{DjMt}, \cite{Vld}. Here, by $\langle
f,\cdot\rangle$, $f\in \mathfrak{D}'(\mathbb{R})$, we denote
sesquilinear functionals on the space $C_{\operatorname
  {comp}}^{\infty}(\mathbb{R})$.

\begin{remark}\label{rm_14_MnRs}
  Without loss of generality, everywhere in the sequel we will
  assume that
\begin{equation*}
  \widehat{q}(0)=0.
\end{equation*}
\end{remark}

Then, the Hill-Schr\"{o}dinger operators can be well defined on the
Hilbert space $L_{2}(\mathbb{R})$ as quasi-differential operators
\cite{SvSh1, SvSh2} by means of the quasi-expressions
\begin{equation*}
  l_{Q}[u]=-(u'-Qu)'-Q(u'-Qu)-Q^{2}u.
\end{equation*}

Set
\begin{align*}
  S_{\operatorname{max}}(q)u & :=l_{Q}[u], \\
  \mathrm{Dom}(S_{\operatorname{max}}(q)) & :=\left\{u\in L_{2}(\mathbb{R})\left|\right.
  u,u'-Qu\in W_{1,{\operatorname {loc}}}^{1}(\mathbb{R}),\,
  l_{Q}[u]\in L_{2}(\mathbb{R})\right\},
\end{align*}
and
\begin{align*}
  \dot{S}_{\operatorname{min}}(q)u & :=l_{Q}[u], \\
  \mathrm{Dom}(\dot{S}_{\operatorname{min}}(q)) & :=\left\{u\in \mathrm{Dom}(S_{\operatorname{max}}(q))\left|\right.
  \mathrm{supp}\,u\Subset\mathbb{R}\right\}.
\end{align*}

It is obvious that the operators $S_{\operatorname{max}}(q)$ are defined
on maximal linear manifolds where the
quasi-expressions $l_{Q}[\cdot]$ are well defined.

\begin{proposition}\label{pr_22_MnRs}
  Let $q(x)\in
  H_{\operatorname{per}}^{-1}\left(\mathbb{R},\mathbb{R}\right)$.
  Then the following statements hold true.
\begin{itemize}
  \item
    [(\rm I)] The operators $\dot{S}_{\operatorname{min}}(q)$ are
    symmetric and lower semibounded on the Hilbert space
    $L_{2}(\mathbb{R})$. In particular, they are closable.
  \item
    [(\rm {II})] The closures $S_{\operatorname{min}}(q)$ of the
    operators $\dot{S}_{\operatorname{min}}(q)$, $
    S_{\operatorname{min}}(q):=(\dot{S}_{\operatorname{min}}(q))^{\sim}$,
    are symmetric, lower semibounded operators on the Hilbert space
    $L_{2}(\mathbb{R})$ with deficiency indices of the form $(m,m)$
    where $0\leq m\leq2$.  The operators $S_{\operatorname{max}}(q)$
    are adjoint to the operators $S_{\operatorname{min}}(q)$,
  \begin{equation*}
    S_{\operatorname{min}}^{\ast}(q)=S_{\operatorname{max}}(q).
  \end{equation*}
  In particular, $S_{\operatorname{max}}(q)$ are closed operators on the Hilbert space
  $L_{2}(\mathbb{R})$, and
  \begin{equation*}
    S_{\operatorname{max}}^{\ast}(q)=S_{\operatorname{min}}(q).
  \end{equation*}
  \item
    [(\rm {III})] Domains $\mathrm{Dom}(S_{\operatorname{min}}(q))$
    of the operators $S_{\operatorname{min}}(q)$ consist of those
    and only those functions $u\in
    \mathrm{Dom}(S_{\operatorname{max}}(q))$ which satisfy the
    conditions
  \begin{equation*}
    [u,v]_{+\infty}-[u,v]_{-\infty}=0\quad \forall v\in
    \mathrm{Dom}(S_{\operatorname{max}}(q)),
  \end{equation*}
  where the limits
  \begin{equation*}
     [u,v]_{+\infty}:=\lim_{x\rightarrow +\infty}[u,v]_{x}
     \quad\text{and}\quad
     [u,v]_{-\infty}:=\lim_{x\rightarrow -\infty}[u,v]_{x}
  \end{equation*}
  are well defined and exist.
\end{itemize}
\end{proposition}

Proposition \ref{pr_22_MnRs}, which describes properties of the
operators $\dot{S}_{\operatorname{min}}(q)$ and
$S_{\operatorname{max}}(q)$, is proved in Appendix A.2 by using
methods of the theory of linear quasi-differential operators.

In Proposition \ref{pr_22.1_MnRs} we define Friedrichs extensions of
the minimal operators $S_{\operatorname{min}}(q)$. But for
convenience we first recall some related facts and prove useful
Lemma~\ref{lm_16_Prl}.

Let $H$ be a Hilbert space, and $\dot{A}$ be a densely defined,
lower semibounded linear operator on $H$. Hence, $\dot{A}$ is a
closable, symmetric operator. Define by $A$ its closure,
$A:=(\dot{A})^{\sim}$.

Set
\begin{equation*}
  \dot{t}[u,v]:=(\dot{A}u,v),\quad
  \mathrm{Dom}(\dot{t}):=\mathrm{Dom}(\dot{A}).
\end{equation*}
As known \cite{Kt}, the sesquilinear form $\dot{t}[u,v]$ is closable,
lower semibounded and symmetric on the Hilbert space $H$. Let
$t[u,v]$ be its closure, $t:=(\dot{t})^{\sim}$.

For the operator $\dot{A}$ there is a uniquely defined its
Friedrichs extension $A_{F}$ \cite{Kt},
\begin{equation*}
  t[u,v]=(A_{F}u,v),\quad
  u\in \mathrm{Dom}(A_{F})\subset\mathrm{Dom}(t),\quad v\in
  \mathrm{Dom}(t).
\end{equation*}
Due to the First Representation Theorem \cite{Kt}, the operator
$A_{F}$ is lower semibounded and self-adjoint. In Lemma
\ref{lm_16_Prl} we describe its domain, but at first note that the
following inclusions take place:
\begin{equation*}
  \dot{A}\subset A\subset A_{F}\subset A^{\ast}.
\end{equation*}

\begin{lemma}\label{lm_16_Prl}
  Let $A_{F}$ be a Friedrichs extension of a densely defined, lower
  semibounded operator $\dot{A}$ on a Hilbert space $H$, and let
  $t[u,v]$ be the densely defined, closed, symmetric, and bounded
  from below sesquilinear form on $H$ constructed from the operator
  $\dot{A}$.  Then
\begin{equation*}
  \mathrm{Dom}(A_{F})=\mathrm{Dom}(t)\cap \mathrm{Dom}(A^{\ast}).
\end{equation*}
\end{lemma}
\begin{proof}
It is obvious that
\begin{equation*}
  \mathrm{Dom}(A_{F})\subset\mathrm{Dom}(t)\cap \mathrm{Dom}(A^{\ast}).
\end{equation*}
Let us prove the inverse inclusion.

Let $u\in \mathrm{Dom}(t)\cap \mathrm{Dom}(A^{\ast})$, and $v\in
\mathrm{Dom}(\dot{A})\subset\mathrm{Dom}(A_{F})\subset\mathrm{Dom}(t)$.
Remark that $\mathrm{Dom}(\dot{A})$ is a core of the form $t[u,v]$
and that $\mathrm{Dom}(t)\cap \mathrm{Dom}(A^{\ast})$ contains
$\mathrm{Dom}(\dot{A})$. Then we have
\begin{equation*}
  (A^{\ast}u,v)=(u,\dot{A}v)=(u,A_{F}v)=\overline{(A_{F}v,u)}=\overline{t[v,u]}=t[u,v],
\end{equation*}
i.e.,
\begin{equation*}
  t[u,v]=(A^{\ast}u,v),\quad u\in \mathrm{Dom}(t)\cap
  \mathrm{Dom}(A^{\ast}),\quad v\in
\mathrm{Dom}(\dot{A}).
\end{equation*}
Due to the First Representation Theorem \cite{Kt} we get that $u\in
\mathrm{Dom}(A_{F})$, i.e.,
\begin{equation*}
  \mathrm{Dom}(t)\cap \mathrm{Dom}(A^{\ast})\subset
  \mathrm{Dom}(A_{F}).
\end{equation*}

The proof is complete.
\end{proof}

\begin{proposition}\label{pr_22.1_MnRs}
  Friedrichs extensions $S_{F}(q)$ of the operators
  $S_{\operatorname{min}}(q)$ are defined in the following way:
  \begin{align*}
    S_{F}(q)u & :=l_{Q}[u], \\
    \mathrm{Dom}(S_{F}(q)) & :=\left\{u\in H^{1}(\mathbb{R})\left|\right.
    u'-Qu\in W_{1,{\operatorname {loc}}}^{1}(\mathbb{R}),\,
    l_{Q}[u]\in L_{2}(\mathbb{R})\right\}.
  \end{align*}
\end{proposition}
\begin{proof}
Let us introduce the sesquilinear forms
\begin{equation*}
  \dot{t}[u,v]:=(\dot{S}_{\operatorname{min}}(q)u,v),\quad
  \mathrm{Dom}(\dot{t}):=\mathrm{Dom}(\dot{S}_{\operatorname{min}}(q)).
\end{equation*}
As is well known \cite{Kt}, the sesquilinear forms $\dot{t}[u,v]$
are densely defined, closable, symmetric and bounded from below on
the Hilbert space $L_{2}(\mathbb{R})$. Taking into account that
$\mathrm{Dom}(\dot{S}_{\operatorname{min}}(q))\subset
H_{\operatorname {comp}}^{1}(\mathbb{R})$, the forms $\dot{t}[u,v]$
can be written as
\begin{equation*}
  \dot{t}[u,v]=(u',v')-(Qu,v')-(Qu',v),\quad
  \mathrm{Dom}(\dot{t})\subset H_{\operatorname {comp}}^{1}(\mathbb{R}).
\end{equation*}
Set
\begin{align*}
  \dot{t}_{1}[u,v] & :=(u',v')+(u,v),\quad
  &\mathrm{Dom}(\dot{t}_{1}):=\mathrm{Dom}(\dot{S}_{\operatorname{min}}(q))\subset H_{\operatorname {comp}}^{1}(\mathbb{R}), \\
  \dot{t}_{2}[u,v] & :=-(Qu,v')-(Qu',v)-(u,v), \quad
  &\mathrm{Dom}(\dot{t}_{2}):=\mathrm{Dom}(\dot{S}_{\operatorname{min}}(q))\subset H_{\operatorname {comp}}^{1}(\mathbb{R}),
\end{align*}
i.e.,
\begin{equation*}
  \dot{t}=\dot{t}_{1}+\dot{t}_{2}.
\end{equation*}
It is well known that the form $\dot{t}_{1}[u,v]$ is closable, and
its closure, $t_{1}[u,v]$, $t_{1}:=(\dot{t}_{1})^{\sim}$, has the
representation
\begin{equation*}
  t_{1}[u,v]=(u',v')+(u,v),\quad
  \mathrm{Dom}(t_{1})=H^{1}(\mathbb{R}).
\end{equation*}
As was shown in \cite{HrMk}, the forms $\dot{t}_{2}[u,v]$ are
$t_{1}$-bounded with relative boundary $0$. So, we finally obtain
that the forms $\dot{t}[u,v]$, which are closures of $t[u,v]$,
$t:=(\dot{t})^{\sim}$, are defined as follows:
\begin{equation*}
  t[u,v]=(u',v')-(Qu,v')-(Qu',v),\quad
  \mathrm{Dom}(t)=H^{1}(\mathbb{R}).
\end{equation*}
And the sesquilinear forms $t[u,v]$ are densely defined, closed,
symmetric, and lower semibounded on the Hilbert space
$L_{2}(\mathbb{R})$.

Further, since
\begin{align*}
  S_{\operatorname{min}}^{\ast}(q)u & =l_{Q}[u], \\
  \mathrm{Dom}(S_{\operatorname{min}}^{\ast}(q))
  & =\left\{u\in L_{2}(\mathbb{R})\left|\right.u,u'-Qu\in W_{1,{\operatorname {loc}}}^{1}(\mathbb{R}),\,
  l_{Q}[u]\in L_{2}(\mathbb{R})\right\},
\end{align*}
applying Lemma \ref{lm_16_Prl} we get the needed representations for
Friedrichs extensions of the operators $\dot{S}_{\operatorname{min}}(q)$.

The proof is complete.
\end{proof}

\begin{statement}\label{st_22.2_MnRs}
  The following inclusions take place:
\begin{equation*}
  \dot{S}_{\operatorname{min}}(q)\subset S_{\operatorname{min}}(q) \subset S_{F}(q) \subset S_{\operatorname{max}}(q)
\end{equation*}
and
\begin{align*}
  \mathrm{Dom}(\dot{S}_{\operatorname{min}}(q)) & \subset H_{\operatorname {comp}}^{1}(\mathbb{R}), \\
  \mathrm{Dom}(S_{\operatorname{min}}(q))  \subset H^{1}(\mathbb{R}), & \quad
  \mathrm{Dom}(S_{F}(q))  \subset H^{1}(\mathbb{R}), \\
  \mathrm{Dom}(S_{\operatorname{max}}(q)) & \subset L_{2}(\mathbb{R})\cap
  H_{loc}^{1}(\mathbb{R}).
\end{align*}
\end{statement}

Statement \ref{st_22.2_MnRs} immediately follows from the
corresponding definitions and not very complicated computations.

Now, our aim is to prove that the maximal quasi-differential
operators $S_{\operatorname{max}}(q)$ are self-adjoint.
\begin{proposition}\label{cr_30_MnRs}
  Let $q(x)\in
  H_{\operatorname{per}}^{-1}\left(\mathbb{R},\mathbb{R}\right)$.
  The following statements are equivalent.
  \begin{itemize}
  \item [(a)] The operators $S_{\operatorname{max}}(q)$ are
    self-adjoint.
  \item [(b)] $\mathrm{Dom}(S_{\operatorname{max}}(q))\subset
    H^{1}(\mathbb{R})$.
  \item [(c)] $u'-Qu\in L_{2}(\mathbb{R})\cap W_{1,{\operatorname
        {loc}}}^{1}(\mathbb{R})\quad \forall u\in
    \mathrm{Dom}(S_{\operatorname{max}}(q))$.
\end{itemize}
\end{proposition}

\begin{proof}
  $(a)$ Let $S_{\operatorname{max}}(q)$ be self-adjoint. Then it
  follows from Proposition \ref{pr_22_MnRs}.II and Statement
  \ref{st_22.2_MnRs} that
  \begin{align*}
    & S_{\operatorname{min}}(q)=S_{F}(q)=S_{\operatorname{max}}(q), \\
    &
    \mathrm{Dom}(S_{\operatorname{min}}(q))=\mathrm{Dom}(S_{F}(q))=\mathrm{Dom}(S_{\operatorname{max}}(q))\subset
    H^{1}(\mathbb{R}),
  \end{align*}
  and $(b)$ is true.
  
  Further, under the assumptions $Q\in
  L_{2,{\operatorname{per}}}(\mathbb{R})$ and $u\in
  H^{1}(\mathbb{R})$ we get that $Qu\in L_{2}(\mathbb{R})$
  \cite{HrMk}, which yields $(c)$.
  
  $(b)$ Let us now assume that
  $\mathrm{Dom}(S_{\operatorname{max}}(q))\subset
  H^{1}(\mathbb{R})$. As above, we get $Qu\in L_{2}(\mathbb{R})$,
  and, as a consequence, $(c)$ follows.  Then statement $(a)$
  follows from the Lagrange identity \eqref{eq_28_Prl}, taking into
  account that
  \begin{equation*}
    [u,v]_{+\infty}=0\quad\text{and}\quad [u,v]_{-\infty}=0
  \end{equation*}
  for $u,v\in L_{2}(\mathbb{R})$ and $u'-Qu,v'-Qv\in
  L_{2}(\mathbb{R})\cap W_{1,{\operatorname
      {loc}}}^{1}(\mathbb{R})$.
  
  $(c)$ Assume that $u'-Qu\in
  L_{2}(\mathbb{R})\cap W_{1,{\operatorname
      {loc}}}^{1}(\mathbb{R})\quad \forall u\in
  \mathrm{Dom}(S_{\operatorname{max}}(q))$. Then applying the
  Lagrange identity \eqref{eq_28_Prl} as above we get $(a)$ and, as
  a consequence, $(b)$.
  
  The proof is complete.
\end{proof}

Hryniv and Mykytyuk \cite{HrMk} studied operators associated via 
the First Representation Theorem \cite{Kt} to the sesquilinear
forms
\begin{equation*}
  t[u,v]=(u',v')-(Qu,v')-(Qu',v),\quad
  \mathrm{Dom}(t)=H^{1}(\mathbb{R}),
\end{equation*}
that is, they have actually studied the Friedrichs extensions
$S_{F}(q)$.

Djakov and Mityagin \cite{DjMt} have also treated the Friedrichs extensions $S_{F}(q)$
\textit{a priori} considering the operators on the domains
\begin{equation*}
  \mathrm{Dom}(S_{F}(q))=\left\{u\in H^{1}(\mathbb{R})\left|\right.
  u'-Qu\in W_{1,{\operatorname {loc}}}^{1}(\mathbb{R}),\, l_{Q}[u]\in
  L_{2}(\mathbb{R})\right\},
\end{equation*}
see Proposition \ref{pr_22.1_MnRs} and Proposition \ref{cr_30_MnRs}.

So, due to Proposition \ref{pr_22_MnRs}.II, we have
\begin{equation*}
  S_{\operatorname{max}}(q)\supset S_{\operatorname{max}}^{\ast}(q),
\end{equation*}
and, therefore, it remains to show that the  operators
$S_{\operatorname{max}}(q)$,
\begin{equation*}
  S_{\operatorname{max}}(q)\subset S_{\operatorname{max}}^{\ast}(q)
\end{equation*}
are symmetric.  We do it by applying Lemma \ref{lm_18_Prl}.

Let us consider the following shift operator on the Hilbert space
$L_{2}(\mathbb{R})$:
\begin{equation*}
  (Uf)(x):=f(x+1),\quad \mathrm{Dom}(U):=L_{2}(\mathbb{R}).
\end{equation*}
 Then
$\sigma_{p}(U)=\varnothing$.

Further, let $f\in \mathrm{Dom}(S_{\operatorname{max}}(q))$. It is obvious that
$Uf\in \mathrm{Dom}(S_{\operatorname{max}}(q))$ too, and it is also true that
\begin{equation*}
  U(S_{\operatorname{max}}(q)f)=Ul_{Q}[f(x)]=l_{Q}[f(x+1)]=l_{Q}[(Uf)(x)]=S_{\operatorname{max}}(q)(Uf),
\end{equation*}
i.e., the operators $S_{\operatorname{max}}(q)$ and $U$ commute.

Taking into account that $S_{\operatorname{max}}(q)$ are the second
order quasi-differential operators, i.e., their possible eigenvalues
cannot have multiplicities more than two, and applying
Lemma~\ref{lm_18_Prl} to the operators $S_{\operatorname{max}}(q)$
and $U$ we obtain the following proposition.

\begin{proposition}\label{pr_24_MnRs}
   The point spectra $\sigma_{p}(S_{\operatorname{max}}(q))$ of the
  quasi-differential operators $S_{\operatorname{max}}(q)$ are
  empty.
\end{proposition}

\begin{theorem}\label{th_26_MnRs}
  The quasi-differential operators $S_{\operatorname{max}}(q)$ are
  self-adjoint.
\end{theorem}

\begin{proof}
  It follows from Proposition \ref{pr_22_MnRs}.II and Proposition
  \ref{pr_24_MnRs} that the minimal symmetric operators
  $S_{\operatorname{min}}(q)$ have deficiency indices of the form
  $(0,0)$, i.e., they are self-adjoint. Due to Proposition
  \ref{pr_22_MnRs}.II, this implies that the operators
  $S_{\operatorname{max}}(q)$ are also self-adjoint.

The proof is complete.
\end{proof}

\begin{corollary}\label{cr_28_MnRs}
  The minimal operators $S_{\operatorname{min}}(q)$, the Friedrichs
  extensions $S_{F}(q)$, and the maximal operators
  $S_{\operatorname{max}}(q)$ coincide. In particular, they are
  self-adjoint and lower semibounded.
\end{corollary}

\begin{corollary}\label{cr_29_MnRs}
  Let $q(x)\in
  H_{\operatorname{per}}^{-1}\left(\mathbb{R},\mathbb{R}\right)$,
  and $q_{n}(x)\in H_{\operatorname{per}}^{-1}\left(\mathbb{R},
    \mathbb{R}\right)$, $n\in \mathbb{N}$, be such that
  \begin{equation*}
    q_{n}(x)\overset{H_{\operatorname{per}}^{-1}\left(\mathbb{R}\right)}{\longrightarrow}
    q(x)\quad\text{as}\quad n\rightarrow\infty.
  \end{equation*}
  Then the Hill-Schr\"{o}dinger operators $S(q_{n})$, $n\in
  \mathbb{N}$, converge to the operators $S(q)$ in the norm
  resolvent sense,
  \begin{equation*}
    \left\|\left(S(q_{n})-\lambda I\right)^{-1}-\left(S(q)-\lambda
        I\right)^{-1}\right\|\rightarrow 0 \quad\text{as}\quad
    n\rightarrow\infty,
  \end{equation*}
  for any $\lambda$ belonging to the resolvent sets of $S(q)$ and
  $S(q_{n})$, $n\in \mathbb{N}$.
\end{corollary}

\begin{proof}
  The proof immediately follows from \cite[Theorem 4.1]{HrMk} and
  Corollary \ref{cr_28_MnRs}.
\end{proof}

In particular, the Hill-Schr\"{o}dinger operators $S(q)$ with
distribution potentials $q(x)\in
H_{\operatorname{per}}^{-1}\left(\mathbb{R}, \mathbb{R}\right)$ are
the limits $S_{\operatorname{lim}}(q)$ of a sequence of operators
$S(q_{n})$, $n\in \mathbb{N}$, with smooth potentials $q_{n}(x)\in
L_{2,{\operatorname{per}}}\left(\mathbb{R}, \mathbb{R}\right)$.  For
instance, taking
\begin{equation*}
  q(x)=\sum_{k\in \mathbb{Z}}\widehat{q}(2k)\,e^{i\,2k\pi x}\in
  H_{\operatorname{per}}^{-1}\left(\mathbb{R},\mathbb{R}\right)
\end{equation*}
one can choose
\begin{equation*}
  q_{n}(x):=\sum_{|k|\leq n}\widehat{q}(2k)\,e^{i\,2k\pi x}\in
  C_{\operatorname{per}}^{\infty}\left(\mathbb{R},\mathbb{R}\right),\quad n\in
  \mathbb{N}.
\end{equation*}

Now, we are going to define the Hill-Schr\"{o}dinger operators with
distribution potentials as form-sum operators \cite{Krt}. We will
show that this definition coincides with the definitions given
above.

Let us consider the following sesquilinear forms on the Hilbert
space $L_{2}(\mathbb{R})$:
\begin{equation*}
  \tau[u,v]:=\left\langle-\frac{d^{2}}{dx^{2}}u,v\right\rangle_{L_{2}(\mathbb{R})}+\left\langle
  q(x)u,v\right\rangle_{L_{2}(\mathbb{R})},\quad
  \mathrm{Dom}(\tau)=H^{1}(\mathbb{R}),
\end{equation*}
generated by the one-dimensional Schr\"{o}dinger operators with $q(x)\in
H_{\operatorname{per}}^{-1}(\mathbb{R},\mathbb{R})$.

Here, $\langle \cdot,\cdot\rangle_{L_{2}(\mathbb{R})}$ denotes the
sesquilinear form on the space $L_{2}(\mathbb{R})$, the spaces
$H^{s}(\mathbb{R})$ and $H^{-s}(\mathbb{R})$ for $s\in \mathbb{R}$,
respectively, which is a (sesquilinear) continuous extension of the
inner product in $L_{2}(\mathbb{R})$~\cite{Brz, Grb},
\begin{equation*}
  \langle
  f,g\rangle_{L_{2}(\mathbb{R})}:=\int_{\mathbb{R}}f(x)\overline{g(x)}\,dx\quad
  \forall f,g\in L_{2}(\mathbb{R}).
\end{equation*}
As is known \cite{Krt}, the sesquilinear forms $\tau[u,v]$ are
densely defined, closed, bounded from below, and are defined on the
Hilbert space $L_{2}(\mathbb{R})$. Due to the First Representation
Theorem \cite{Kt}, there are associated operators
$S_{\operatorname{form}}(q)$ that are uniquely defined on the
Hilbert space $L_{2}(\mathbb{R})$, self-adjoint, lower semibounded,
and such that
\begin{itemize}
  \item [i)] $\mathrm{Dom}\left(S_{\operatorname{form}}(q)\right)\subset
  \mathrm{Dom}\left(\tau\right)$ and
  \begin{equation*}
     \tau[u,v]=(S_{\operatorname{form}}(q)u,v)\quad \forall u\in
     \mathrm{Dom}\left(S_{\operatorname{form}}(q)\right),\quad \forall v\in
     \mathrm{Dom}\left(\tau\right);
  \end{equation*}
  \item [ii)] $\mathrm{Dom}\left(S_{\operatorname{form}}(q)\right)$ are cores of
  the forms $\tau[u,v]$;
  \item [iii)] if $u\in \mathrm{Dom}\left(\tau\right)$, $w\in
  L_{2}(\mathbb{R})$, and
  \begin{equation*}
    \tau[u,v]=(w,v)
  \end{equation*}
  holds for every $v$ in cores of the forms $\tau[u,v]$,
  then $u\in \mathrm{Dom}\left(S_{\operatorname{form}}(q)\right)$ and
  \begin{equation*}
    S_{\operatorname{form}}(q)u=w.
  \end{equation*}
\end{itemize}
The operators $S_{\operatorname{form}}(q)$ are called form-sum
operators associated with the forms $\tau[u,v]$, and denoted by
\begin{equation*}
  S_{\operatorname{form}}(q):=-\frac{d^{2}}{dx^{2}}\dotplus q(x).
\end{equation*}
It will also be convenient to use the notations
\begin{equation*}
  \tau_{S_{\operatorname{form}}(q)}[u,v]\equiv \tau[u,v].
\end{equation*}

\begin{proposition}{\rm{(\cite{Krt})}}\label{pr_32_MnRs}.
The Hill-Schr\"{o}dinger operators with  distribution potentials
from the negative Sobolev space $H_{\operatorname{per}}^{-1}(\mathbb{R},
\mathbb{R})$ are well defined on the Hilbert space
$L_{2}(\mathbb{R})$ as self-adjoint, lower semibounded form-sum
operators $S_{\operatorname{form}}(q)$,
\begin{equation*}
  S_{\operatorname{form}}(q)=-\frac{d^{2}}{dx^{2}}\dotplus q(x),
\end{equation*}
associated with the sesquilinear forms
\begin{equation*}
  \tau_{S_{\operatorname{form}}(q)}[u,v]=\bigg\langle-\frac{d^{2}}{dx^{2}}u,v\bigg\rangle_{L_{2}(\mathbb{R})}+\left\langle
  q(x)u,v\right\rangle_{L_{2}(\mathbb{R})},\quad
  \mathrm{Dom}(\tau)=H^{1}(\mathbb{R}),
\end{equation*}
acting on the dense domains
\begin{equation*}
  \mathrm{Dom}\left(S_{\operatorname{form}}(q)\right):=\left\{u\in H^{1}(\mathbb{R})\left|
  -\frac{d^{2}}{dx^{2}}u+q(x)u\in L_{2}(\mathbb{R})\right.\right\}
\end{equation*}
as
\begin{equation*}
  S_{\operatorname{form}}(q)u:=-\frac{d^{2}}{dx^{2}}u+q(x)u\in L_{2}(\mathbb{R}),\quad u\in
  \mathrm{Dom}\left(S_{\operatorname{form}}(q)\right).
\end{equation*}
\end{proposition}

\begin{theorem}\label{th_34_MnRs}
  The quasi-differential operators $S(q)$ and the form-sum operators
  $S_{\operatorname{form}}(q)$ coincide.
\end{theorem}

\begin{proof}
  Let $u\in \mathrm{Dom}\left(S(q)\right)$. Recall that
\begin{equation*}
  \mathrm{Dom}(S(q))=\left\{u\in H^{1}(\mathbb{R})\left|\right.u'-Qu\in W_{1,{\operatorname {loc}}}^{1}(\mathbb{R}),\,
  l_{Q}[u]\in L_{2}(\mathbb{R})\right\},
\end{equation*}
so that
\begin{equation*}
  \mathrm{Dom}(S(q))\subset
  \mathrm{Dom}(\tau_{S_{\operatorname{form}}(q)})=H^{1}(\mathbb{R}).
\end{equation*}
Then we have
\begin{align*}
  \tau_{S_{\operatorname{form}}(q)}[u,v] & =\langle -u'',v\rangle_{L_{2}(\mathbb{R})}+\langle
  q(x)u,v\rangle_{L_{2}(\mathbb{R})}
  = \langle u',v'\rangle_{L_{2}(\mathbb{R})}-\langle Q(x),\overline{u'}v+\overline{u}v'\rangle_{L_{2}(\mathbb{R})} \\
  & = (u',v')-(Qu,v')-(Qu',v)=(l_{Q}[u],v)\quad \forall v\in
  C_{\operatorname {comp}}^{\infty}(\mathbb{R}).
\end{align*}
And, due to the First Representation Theorem \cite{Kt}, we conclude
that
\begin{equation*}
  u\in \mathrm{Dom}(S_{\operatorname{form}}(q))\quad\text{and}\quad S_{\operatorname{form}}(q)u=l_{Q}[u],
\end{equation*}
i.e.,
\begin{equation*}
  S(q)\subset S_{\operatorname{form}}(q).
\end{equation*}
Taking into account that the operators $S(q)$ and
$S_{\operatorname{form}}(q)$ are self-adjoint, the latter also gives
the inverse inclusions
\begin{equation*}
  S(q)\supset S_{\operatorname{form}}(q).
\end{equation*}
The proof is complete.
\end{proof}

\subsection{Spectra of the Hill-Schr\"{o}dinger operators with
distribution potentials} 

In this section, we will establish characteristic properties of the
structure of the spectrum of the Hill-Schr\"{o}dinger operators
$S(q)$ with distribution potentials $q(x)\in
H_{\operatorname{per}}^{-1}\left(\mathbb{R},\mathbb{R}\right)$.
Using a limit process in the generalized sense applied to the
Hill-Schr\"{o}dinger operators $S(q_{n})$, $n\in \mathbb{N}$, with
smooth potentials $q_{n}(x)\in
L_{2,{\operatorname{per}}}(\mathbb{R},\mathbb{R})$ (see Corollary
\ref{cr_29_MnRs}) we show that the Hill-Schr\"{o}dinger operators
$S(q)$, with the distributions $q(x)\in H_{\operatorname{per}}^{-1}
\left(\mathbb{R}, \mathbb{R}\right)$ as potentials, have continuous
spectra with a band and gap structure.

For different approaches, see \cite{HrMk, Krt, DjMt}.

At first, let us recall well known results related to the
classical case of
$L_{2,{\operatorname{per}}}(\mathbb{R},\mathbb{R})$-potentials
$q(x)$,
\begin{equation}\label{eq_38_MnRs}
  q(x)\in L_{2,{\operatorname{per}}}(\mathbb{R},\mathbb{R}),
\end{equation}
see, for an example, \cite{DnSch2, ReSi4}. Under assumption
\eqref{eq_38_MnRs}, the Hill-Schr\"{o}dinger operators $S(q)$ are
lower semibounded and self-adjoint on the Hilbert space
$L_{2}(\mathbb{R})$; they have absolutely continuous spectra with a
band and gap structure.

Spectra of the Hill-Schr\"{o}dinger operators are well defined by
locating the spectrum gap endpoints. It is known that for the endpoints
$\{\lambda_{0}(q),\lambda_{k}^{\pm}(q)\}_{k=1}^{\infty}$ of the
spectrum gaps, we have the following inequalities:
\begin{equation}\label{eq_40_MnRs}
  -\infty<\lambda_{0}(q)<\lambda_{1}^{-}(q)\leq\lambda_{1}^{+}(q)<
  \lambda_{2}^{-}(q)\leq\lambda_{2}^{+}(q)<\cdots\,
\end{equation}
The spectrum bands (or stability zones),
\begin{equation*}
  \mathcal{B}_{0}(q):=[\lambda_{0}(q),\lambda_{1}^{-}(q)],\quad
  \mathcal{B}_{k}(q):=[\lambda_{k}^{+}(q),\lambda_{k+1}^{-}(q)],\quad k\in
  \mathbb{N},
\end{equation*}
are characterized as a set of real $\lambda\in \mathbb{R}$
for which all solutions of the equation
\begin{equation}\label{eq_42_MnRs}
  S(q) u=\lambda u
\end{equation}
are bounded. On the other hand, spectrum gaps (or instability
zones),
\begin{equation*}
  \mathcal{G}_{0}(q):=(-\infty,\lambda_{0}(q)),\quad
  \mathcal{G}_{k}(q):=(\lambda_{k}^{-}(q),\lambda_{k}^{+}(q)),\quad k\in
  \mathbb{N},
\end{equation*}
make a set of real $\lambda\in \mathbb{R}$ for which any nontrivial
solution of the equation \eqref{eq_42_MnRs} is unbounded.

As follows from~\eqref{eq_40_MnRs},  it could happen that
\begin{equation*}
  \lambda_{k}^{-}(q)=\lambda_{k}^{+}(q)
\end{equation*}
for some $k\in \mathbb{N}$. In such a case, we say that the
corresponding spectrum gap $\mathcal{G}_{k}(q)$ is
\textit{collapsed} or \textit{closed}. Note that this cannot happen
for spectrum bands.

Further, it could happen that the endpoints of spectrum gaps for
even numbers $k\in \mathbb{Z}_{+}$ are periodic eigenvalues of the
problem on the interval $[0,1]$,
\begin{equation*}
   S_{+}(q)u:=-u''+q(x)u=\lambda u,\quad u\in
   \mathrm{Dom}\left(S_{+}(q)\right),
\end{equation*}
and the endpoints of spectrum gaps for odd numbers $k\in \mathbb{N}$ are
semiperiodic eigenvalues of the problem on the interval $[0,1]$,
\begin{equation*}
   S_{-}(q)u:=-u''+q(x)u=\lambda u,\quad u\in
   \mathrm{Dom}\left(S_{-}(q)\right).
\end{equation*}
Under the assumption \eqref{eq_38_MnRs}, domains of the operators $S_{+}(q)$ and
$S_{-}(q)$ have the form
\begin{equation*}
  \mathrm{Dom}(S_{\pm}(q))=\left\{u\in H^{2}[0,1]\left|\,
  u^{(j)}(0)=\pm u^{(j)}(1),\, j=0,1\right.\right\}.
\end{equation*}

Now, applying the limit process in the generalized sense (see Corollary
\ref{cr_29_MnRs}) to the Hill-Schr\"{o}dinger operators $S(q_{n})$, $n\in
\mathbb{N}$, with $L_{2,{\operatorname{per}}}(\mathbb{R}, \mathbb{R})$-potentials $q_{n}(x)$ we
establish the following statement.

\begin{theorem}\label{th_36_MnRs}
  Suppose that $q(x)\in H_{\operatorname{per}}^{-1}\left(\mathbb{R},
    \mathbb{R}\right)$. Then the Hill-Schr\"{o}dinger operators
  $S(q)$ have continuous spectra with a band and gap structure such
  that the endpoints $\{\lambda_{0}(q),
  \lambda_{k}^{\pm}(q)\}_{k=1}^{\infty}$ of the spectrum gaps satisfy
  the inequalities
\begin{equation*}
  -\infty<\lambda_{0}(q)<\lambda_{1}^{-}(q)\leq\lambda_{1}^{+}(q)<
  \lambda_{2}^{-}(q)\leq\lambda_{2}^{+}(q)<\cdots\,
\end{equation*}
Moreover, the endpoints of the spectrum gaps for even (odd) numbers $k\in \mathbb{Z}_{+}$
are periodic (semiperiodic) eigenvalues of the problem on the interval $[0,1]$,
\begin{equation*}
   S_{\pm}(q)u=-u''+q(x)u=\lambda u,\quad u\in
   \mathrm{Dom}\left(S_{\pm}(q)\right).
\end{equation*}
\end{theorem}
\begin{remark}
  The operators $S_{+}(q)$ and $S_{-}(q)$ are well defined on the
  Hilbert space $L_{2}(0,1)$ as lower semi-bounded, self-adjoint
  form-sum operators,
\begin{equation*}
  S_{\pm}(q)=\left(-\frac{d^{2}}{dx^{2}}\right)_{\pm}\dotplus q(x).
\end{equation*}
They also can be well defined in alternative equivalent ways, --- as
quasi-differential operators or as limits, in the norm resolvent
sense, of a sequence of operators with smooth potentials.

In the papers \cite{MkhMl3, MkhMl4, MkhMl5}, the authors
meticulously treated the form-sum operators
\begin{equation*}
  S_{\pm}(V)=\left((-1)^{m}\frac{d^{2m}}{dx^{2m}}\right)_{\pm}\dotplus
  V(x),\quad V(x)\in H_{\operatorname{per}}^{-m}[0,1],\quad m\in
  \mathbb{N},
\end{equation*}
defined on $L_{2}(0,1)$.

In \cite{Mlb, MkhMl1, MkhMl2}, the authors studied two terms
differential operators of even order defined in the
\textit{negative} Sobolev spaces.
\end{remark}

\begin{proof}
  Let $\left\{q_{n}(x)\right\}_{n\in\mathbb{N}}$ be a sequence of
  real-valued trigonometric polynomials, which converges to the
  singular potential $q(x)$ in the norm of the space
  $H_{\operatorname{per}}^{-1}\left(\mathbb{R}\right)$. With this
  sequence one can associate a sequence of self-adjoint operators
  $\left\{S_{\pm}(q_{n})\right\}_{n\in \mathbb{N}}$ defined in
  $L_{2}(0,1)$, and a sequence of Hill operators
  $\left\{S(q_{n})\right\}_{n\in \mathbb{N}}$ defined on
  $L_{2}(\mathbb{R})$. As was proved by the authors in \cite{MkhMl3,
    MkhMl5}, the sequences $\left\{S_{\pm}(q_{n})\right\}_{n\in
    \mathbb{N}}$ converge to the operators $S_{\pm}(q)$ in the norm
  resolvent sense. Hence, eigenvalues of these operators
  $\left\{S_{\pm}(q_{n})\right\}_{n\in \mathbb{N}}$ converge to the
  corresponding eigenvalues of the limit operators $S_{\pm}(q)$
  \cite[Theorem VIII.23 and Theorem~VIII.24]{ReSi1} (also see
  \cite{Kt}). Further, as is well known \cite{CdLv, DnSch2}, for the
  operators $\left\{S_{\pm}(q_{n})\right\}_{n\in \mathbb{N}}$, the
  assertion of the theorem is true, i.e.,
\begin{equation}\label{eq_44_MnRs}
  -\infty<\lambda_{0}(q_{n})<\lambda_{1}^{-}(q_{n})\leq\lambda_{1}^{+}(q_{n})<
  \lambda_{2}^{-}(q_{n})\leq\lambda_{2}^{+}(q_{n})<\cdots\,
\end{equation}
Moreover, as we have already proved (see Corollary
\ref{cr_29_MnRs}), the sequence $\left\{S(q_{n})\right\}_{n\in
  \mathbb{N}}$ converges to the operator $S(q)$ in the norm resolvent
sense. Therefore, from \eqref{eq_44_MnRs} we get
\begin{equation*}
  -\infty<\lambda_{0}(q)\leq \lambda_{1}^{-}(q)\leq\lambda_{1}^{+}(q)\leq
  \lambda_{2}^{-}(q)\leq\lambda_{2}^{+}(q)\leq\cdots\, ,
\end{equation*}
where $\lambda_{0}(q),\,\lambda_{2k}^{\pm}(q)\in \sigma(S_{+}(q))$
and $\lambda_{2k-1}^{\pm}(q)\in \sigma(S_{-}(q))$, $k\in
\mathbb{N}$.

Now it remains to show that the strict inequalities
\begin{equation*}
    \lambda_{k}^{+}(q_{n})< \lambda_{k+1}^{-}(q_{n}),\quad k\in
    \mathbb{Z}_{+},
\end{equation*}
can not become equalities. Indeed, suppose the contrary. Then, one
of the spectrum zones of the operator $S(q)$ degenerates into a
point,
\begin{equation*}
    \lambda_{k_{0}}^{+}(q)=\lambda_{k_{0}+1}^{-}(q),\quad k_{0}\in
    \mathbb{Z}_{+}.
\end{equation*}
Since it is an isolated point of the spectrum of the operator
$S(q)$, it cannot belong to the continuous spectrum
$\sigma_{c}\left(S(q)\right)$. On the other hand, it cannot belong
to the point spectrum of the operator $S(q)$, since
$\sigma_{p}\left(S(q)\right)=\varnothing$. The obtained
contradiction proves the inequalities in theorem.

The proof is complete.
\end{proof}

\section{Concluding remarks}\label{sec_CnRm}

It follows from the direct integral decomposition of the
Hill-Schr\"{o}dinger operators $S(q)$ \cite{HrMk} and \cite[Theorem
XIII.86]{ReSi4} that $\sigma_{sc}(S(q))=\varnothing$.  Therefore,
the continuity of spectra of the operators $S(q)$, which was proved
in this paper, shows that they are absolutely continuous \cite{MiSb}.

From Theorem C and the results of the authors in~\cite{MkhMl3}, one
obtains a series of theorems establishing relationships between the
lengths of the spectrum gaps and smoothness of the distribution
potentials $q(x)\in
H_{\operatorname{per}}^{-s}(\mathbb{R},\mathbb{R})$, $s\geq -1$, of
the Hill-Schr\"{o}dinger operators $S(q)$ \cite{MkhMl6}.

\bigskip {\it{Acknowledgments.}} The investigation of the first
author was partially supported by the Ukrainian Foundation for
Fundamental Research, Grant 14.1/003.

\section*{Appendix: Some proofs}\label{sec_SmPr_A}
\noindent A.1. \textbf{Proof of Statement \ref{st_15.1_Prl}.}\label{prf_A10} 
At first note that  the relations
\begin{equation*}
  \dot{L}_{\operatorname{min}}(q)\subset L_{\operatorname{min}}(q)
\end{equation*}
give
\begin{equation*}
  (\dot{L}_{\operatorname{min}}(q))^{\sim}\subset L_{\operatorname{min}}(q),
\end{equation*}
see Proposition \ref{pr_15_Prl}.III. Therefore, it suffices to show the inverse
inclusions,
\begin{equation*}
 (\dot{L}_{\operatorname{min}}(q))^{\sim}\supset L_{\operatorname{min}}(q).
\end{equation*}

Let $\Delta=[\alpha,\beta]$ denote a fixed, closed interval that
completely lies in the interval $[0,1]$, and let
\begin{equation*}
  \mathfrak{H}_{\Delta}:=L_{2}(\alpha,\beta).
\end{equation*}
On the Hilbert space $\mathfrak{H}_{\Delta}$, consider the operators
$L_{{\operatorname{min}},\Delta}(q)$ and
$L_{{\operatorname{max}},\Delta}(q)$ generated by $l_{Q}[\cdot]$ on
the interval $\Delta$, which are are mutually adjoint due to
Proposition \ref{pr_15_Prl}.III,
\begin{equation*}
  L_{{\operatorname{min}},\Delta}^{\ast}(q)=L_{{\operatorname{max}},\Delta}(q),\quad
  L_{{\operatorname{max}},\Delta}^{\ast}(q)=L_{{\operatorname{min}},\Delta}(q).
\end{equation*}

On the other hand the Hilbert space $\mathfrak{H}_{\Delta}$ can be
well embedded into the space $\mathfrak{H}:=L_{2}(0,1)$ assuming
that the function $u\in \mathfrak{H}_{\Delta}$ equals zero on the
interval $\Delta$. Thus, the domains
$\mathrm{Dom}(L_{{\operatorname{min}},\Delta}(q))$ of the operators
$L_{{\operatorname{min}},\Delta}(q)$ become a part of the domains
$\mathrm{Dom}(L_{\operatorname{max}}(q))$ of the operators
$L_{\operatorname{max}}(q)$, since continuity of the
quasi-derivatives $u^{[j]}(x)$, $j=0,1$, of the function $u\in
\mathrm{Dom}(L_{{\operatorname{min}},\Delta}(q))$ is preserved when
extending the function over the interval $\Delta$. Moreover,
extended in such a way, the function $u\in
\mathrm{Dom}(L_{{\operatorname{min}},\Delta}(q))$ then belongs to
$\mathrm{Dom}(\dot{L}_{\operatorname{min}}(q))$. Therefore, if $v\in
\mathrm{Dom}(\dot{L}_{\operatorname{min}}^{\ast}(q))$, then we have
\begin{equation}\label{eq_26.1_Prl}
  \left(\dot{L}_{\operatorname{min}}^{\ast}(q)v,u\right)=\left(v,\dot{L}_{\operatorname{min}}(q)u\right)\quad
  \forall u\in \mathrm{Dom}(L_{{\operatorname{min}},\Delta}(q)).
\end{equation}
Since $u(x)=0$ on the interval $\Delta$, the scalar product in
\eqref{eq_26.1_Prl} is the $\mathfrak{H}_{\Delta}$-inner product.
Denoting these scalar products with the index $\Delta$ we can
rewrite \eqref{eq_26.1_Prl} as follows:
\begin{equation*}
  \left((\dot{L}_{\operatorname{min}}^{\ast}(q)v)_{\Delta},u\right)_{\Delta}=
  \left(v_{\Delta},L_{{\operatorname{min}},\Delta}(q)u\right)_{\Delta}\quad
  \forall u\in \mathrm{Dom}(L_{{\operatorname{min}},\Delta}(q)).
\end{equation*}
Here, $(\dot{L}_{\operatorname{min}}^{\ast}(q)v)_{\Delta}$,
$v_{\Delta}$ denote the functions
$\dot{L}_{\operatorname{min}}^{\ast}(q)v$ and $v$ considered only
in the interval $\Delta$. So, from the latter we obtain
\begin{equation*}
  v_{\Delta}\in \mathrm{Dom}(L_{{\operatorname{min}},\Delta}^{\ast}(q))=\mathrm{Dom}(L_{{\operatorname{max}},\Delta}(q))
\end{equation*}
and
\begin{equation*}
  (\dot{L}_{\operatorname{min}}^{\ast}(q)v)_{\Delta}=L_{{\operatorname{min}},\Delta}^{\ast}(q)v_{\Delta}=L_{{\operatorname{max}},\Delta}(q)v_{\Delta}=\left(l_{Q}[v]\right)_{\Delta}.
\end{equation*}
Since these relations hold for any interval $\Delta\subset [0,1]$, we
conclude that
\begin{equation*}
  v\in \mathrm{Dom}(L_{\operatorname{max}}(q))\quad\text{and}\quad
  \dot{L}_{\operatorname{min}}^{\ast}(q)v=l_{Q}[v]=L_{\operatorname{max}}(q)v.
\end{equation*}
Thus, we have proved that
\begin{equation*}
 \dot{L}_{\operatorname{min}}^{\ast}(q)\subset L_{\operatorname{max}}(q),
\end{equation*}
i.e.,
\begin{equation*}
  \dot{L}_{\operatorname{min}}^{\ast\ast}(q)\supset
  L_{\operatorname{max}}^{\ast}(q)=L_{\operatorname{min}}(q),
\end{equation*}
which implies the required inclusions
\begin{equation*}
  (\dot{L}_{\operatorname{min}}(q))^{\sim}\supset L_{\operatorname{min}}(q).
\end{equation*}

The proof is complete.
$\hfill \square$

\medskip


\noindent A.2. \textbf{Proof of Proposition \ref{pr_22_MnRs}.}\label{prf_A12}
$(\verb"I")$ At first note that
\begin{equation}\label{eq_34_MnRs}
  \mathrm{Dom}(\dot{S}_{\operatorname{min}}(q))\subset
  H_{\operatorname {comp}}^{1}(\mathbb{R}).
\end{equation}
Let $u\in \mathrm{Dom}(\dot{S}_{\operatorname{min}}(q))$. Then we have
\begin{equation*}
  (\dot{S}_{\operatorname{min}}(q)u,u)=(l_{Q}[u],u)=(u',u')-(Qu,u')-(Qu',u),
\end{equation*}
taking into account that, due to the \eqref{eq_34_MnRs},
\begin{equation*}
  |u'|^{2},\, Quu'\in L_{1,{\operatorname{comp}}}(\mathbb{R}).
\end{equation*}
Now, we estimate $(Qu,u')$ and $(Qu',u)$ as in \cite{HrMk},
\begin{equation*}
   \left|(Qu,u')\right| \leq \|Q\|_{L_{2,{\operatorname{per}}}(\mathbb{R})} \left(\varepsilon\|u'\|_{L_{2}(\mathbb{R})}+
   b(\varepsilon^{-1})\|u\|_{L_{2}(\mathbb{R})}\right),\quad \varepsilon\in
   (0,1],\quad b\geq 0,
\end{equation*}
which yields
\begin{equation*}
  (\dot{S}_{\operatorname{min}}(q)u,u)\geq
  -\gamma(\varepsilon^{-1})\|u\|_{L_{2}(\mathbb{R})}\quad
  \forall u\in \mathrm{Dom}(\dot{S}_{\operatorname{min}}(q)),\quad \gamma\geq 0.
\end{equation*}
We can conclude that $\dot{S}_{\operatorname{min}}(q)$ are Hermitian
operators, lower semibounded on $L_{2}(\mathbb{R})$.

Now, let us show that
$\mathrm{Dom}(\dot{S}_{\operatorname{min}}(q))$ are dense in the
Hilbert space $L_{2}(\mathbb{R})$.

Obviously, it is sufficient to prove that any element $h\in
\mathfrak{H}$, $\mathfrak{H}:=L_{2}(\mathbb{R})$, which is
orthogonal to $\mathrm{Dom}(\dot{S}_{\operatorname{min}}(q))$ is
equal to zero.  Suppose that $h(x)$ is  such a function,
\begin{equation*}
  h(x)\perp \mathrm{Dom}(\dot{S}_{\operatorname{min}}(q)),
\end{equation*}
and let $\Delta=[\alpha,\beta]$ be a fixed, closed interval
compactly lying in the real axis $\mathbb{R}$ ($\Delta\Subset
\mathbb{R}$). Any element $u\in
\mathrm{Dom}(S_{{\operatorname{min}},\Delta}(q))$ can be viewed as
an element of $\mathrm{Dom}(\dot{S}_{\operatorname{min}}(q))$ (for
the notations see the proof of Statement \ref{st_15.1_Prl}),
consequently, $h(x)$ is orthogonal to
$\mathrm{Dom}(S_{{\operatorname{min}},\Delta}(q))$. Due to
Proposition \ref{pr_15_Prl}.II,
$\mathrm{Dom}(S_{{\operatorname{min}},\Delta}(q))$ is dense in
$\mathfrak{H}_{\Delta}=L_{2}(\alpha,\beta)$, hence the function
$h(x)$ considered in the interval $\Delta$ has to be equal to zero
almost everywhere in $\Delta$.

Since the interval $\Delta\Subset \mathbb{R}$ was arbitrary, we
conclude that $h(x)=0$ almost everywhere on $\mathbb{R}$.

So, statement $(\verb"I")$ of Proposition \ref{pr_22_MnRs} has been
proved completely.

$(\verb"II")$ It is obvious that the operators
$S_{\operatorname{min}}(q)$ are symmetric, lower semibounded on the
Hilbert space $L_{2}(\mathbb{R})$.

Let us show that the operators $S_{\operatorname{min}}(q)$ and
$S_{\operatorname{max}}(q)$ are adjoint to each other. Since
$(\dot{S}_{\operatorname{min}}(q))^{\sim}=S_{\operatorname{min}}(q)$,
we have
$\dot{S}_{\operatorname{min}}^{\ast}(q)=S_{\operatorname{min}}^{\ast}(q)$,
and it suffices to show that
\begin{equation*}
  \dot{S}_{\operatorname{min}}^{\ast}(q)=S_{\operatorname{max}}(q).
\end{equation*}
Applying the Lagrange identity \eqref{eq_28_Prl}, we have
\begin{equation*}
  (S_{\operatorname{max}}(q)u,v)=(u,\dot{S}_{\operatorname{min}}(q)v)\quad \forall u\in
  \mathrm{Dom}(S_{\operatorname{max}}(q)),\quad  \forall v\in \mathrm{Dom}(\dot{S}_{\operatorname{min}}(q)),
\end{equation*}
which implies that
\begin{equation*}
  S_{\operatorname{max}}(q)\subset \dot{S}_{\operatorname{min}}^{\ast}(q).
\end{equation*}
So, it remains to prove the inverse inclusions,
\begin{equation*}
  S_{\operatorname{max}}(q)\supset \dot{S}_{\operatorname{min}}^{\ast}(q).
\end{equation*}
We do it in a similar manner as in the proof of Statement
\ref{st_15.1_Prl}.

Let $v(x)$ be an arbitrary element in the domains
$\mathrm{Dom}(\dot{S}_{\operatorname{min}}^{\ast}(q))$ of the
operators $\dot{S}_{\operatorname{min}}^{\ast}(q)$, and let
$\Delta=[\alpha,\beta]$ be a fixed, compact interval ($\Delta\Subset
\mathbb{R}$). As in the proof of Statement \ref{st_15.1_Prl}, we
obtain
\begin{equation*}
  \left((\dot{S}_{\operatorname{min}}^{\ast}(q)v)_{\Delta},u\right)_{\Delta}=
  \left(v_{\Delta},S_{{\operatorname{min}},\Delta}(q)u\right)_{\Delta}\quad
  \forall u\in \mathrm{Dom}(S_{{\operatorname{min}},\Delta}(q)).
\end{equation*}
So, one can conclude that
\begin{equation*}
  v_{\Delta}\in \mathrm{Dom}(S_{{\operatorname{max}},\Delta}(q))
\end{equation*}
and
\begin{equation*}
  (\dot{S}_{\operatorname{min}}^{\ast}(q)v)_{\Delta}=S_{{\operatorname{min}},\Delta}^{\ast}(q)v_{\Delta}
  =S_{{\operatorname{max}},\Delta}(q)v_{\Delta}=\left(l_{Q}[v]\right)_{\Delta}.
\end{equation*}
Taking into account that the interval $\Delta\subset \mathbb{R}$ is
arbitrarily chosen, we finally get that
\begin{equation*}
  v\in \mathrm{Dom}(S_{\operatorname{max}}(q))\quad\text{and}\quad
  \dot{S}_{\operatorname{min}}^{\ast}(q)v=l_{Q}[v]=S_{\operatorname{max}}(q)v,
\end{equation*}
so that the required inclusions hold,
\begin{equation*}
  S_{\operatorname{max}}(q)\supset \dot{S}_{\operatorname{min}}^{\ast}(q).
\end{equation*}

Further, let us find the deficiency index of the operators
$S_{\operatorname{min}}(q)$. At first it is necessary to note that,
since the operators $S_{\operatorname{min}}(q)$ are lower
semibounded, their deficiency indices are equal.

Let $\lambda\in \mathbb{C}$, $\mathrm{Im}\,\lambda\neq 0$. Then the
deficiency indices of the operators $S_{\operatorname{min}}(q)$,
which  will be denoted by $m$, are equal to the number of linearly
independent solutions of the equation
\begin{equation*}
  S_{\operatorname{min}}^{\ast}(q)u=\lambda u,
\end{equation*}
i.e., of the equation (Proposition \ref{pr_22_MnRs}.II)
\begin{equation*}
  S_{\operatorname{max}}(q)u=\lambda u.
\end{equation*}
In other words the deficiency index is a maximal number of linear
independent solutions of the equation
\begin{equation*}
  l_{Q}[u]=\lambda u
\end{equation*}
in the Hilbert space $L_{2}(\mathbb{R})$. Since the total number of
linearly independent solutions of this equation is $2$, we conclude
that
\begin{equation*}
  0\leq m\leq 2.
\end{equation*}

Assertion $(\verb"II")$ is proved.

$(\verb"III")$ Let $u,v\in \mathrm{Dom}(S_{\operatorname{max}}(q))$.
Then applying the Lagrange identity \eqref{eq_28_Prl} we conclude
that the following limits exist:
\begin{equation*}
   [u,v]_{+\infty}:=\lim_{x\rightarrow +\infty}[u,v]_{x}
   \quad\text{and}\quad
   [u,v]_{-\infty}:=\lim_{x\rightarrow -\infty}[u,v]_{x},
\end{equation*}
and, as a consequence, the Lagrange identity \eqref{eq_28_Prl} becomes
\begin{equation}\label{eq_36_MnRs}
  (l_{Q}[u],v)-(u,l_{Q}[v])=[u,v]_{-\infty}^{+\infty}\quad
  \forall u,v\in \mathrm{Dom}(S_{\operatorname{max}}(q)).
\end{equation}
Now, due to Proposition \ref{pr_22_MnRs}.II, we have
\begin{equation*}
  S_{\operatorname{min}}(q)=S_{\operatorname{max}}^{\ast}(q).
\end{equation*}
Therefore, the domains $\mathrm{Dom}(S_{\operatorname{min}}(q))$
consist of only the functions $u\in
\mathrm{Dom}(S_{\operatorname{max}}(q))$ that satisfy
the identities
\begin{equation*}
  (u,S_{\operatorname{max}}(q)v)=(S_{\operatorname{max}}(q)u,v)\quad \forall v\in \mathrm{Dom}(S_{\operatorname{max}}(q))
\end{equation*}
and only of them.  Together with the Lagrange identity
\eqref{eq_36_MnRs}, the latter implies the required assertion, i.e.,
\begin{equation*}
  u\in \mathrm{Dom}(S_{\operatorname{min}}(q))\Leftrightarrow
  [u,v]_{+\infty}-[u,v]_{-\infty}=0,\quad u\in
  \mathrm{Dom}(S_{\operatorname{max}}(q))\quad
  \forall v\in \mathrm{Dom}(S_{\operatorname{max}}(q)).
\end{equation*}

Proposition \ref{pr_22_MnRs} is proved.  $\hfill \square $

\end{document}